\documentclass[12pt]{amsart}
\usepackage{latexsym,amsmath,amssymb,amscd}
\begin{document}
\title{Reduced $C^*$-algebra of the $p$-adic group $GL(n)$ II}
\author{R.J. Plymen}
\date{Version 1.5: \today}
\maketitle
\section{Introduction}
Harish-Chandra, in his search for the Plancherel measure on a
reductive $p$-adic group, had a clear conception of the {\it
support} of the Plancherel measure \cite{HC}.   The support is a
$C^{\infty}$-manifold with countably many connected components.
Each component is a compact torus $T$ and Plancherel measure is
absolutely continuous with respect to Haar measure on $T$.  The
Plancherel density is a real analytic function which is invariant
under a certain finite group which acts on $T$, see \cite[p. 353
-- 367]{HC}.

Let $F$ be a nonarchimedean local field, let $G = GL(n) =
GL(n,F)$, and let $C^*_r(G)$ denote the reduced $C^*$-algebra of
$G$. In \cite{LP} \cite{P1} we proved that $C^*_r(G)$ is strongly
Morita equivalent to a commutative $C^*$-algebra:
\[
 C^*_r(G) \sim C_0(\Omega^t(G))
\]
where $\Omega^t(G)$ is the Harish-Chandra parameter space.  In
\cite{P2}, we constructed an equivalence bimodule which effects
this strong Morita equivalence.   For equivalence bimodules and
strong Morita equivalence, see \cite[II.Appendix A, p.152]{Co}.

Let $\mathcal{H}(G)$ be the Hecke algebra of $G$, and let $\Omega$
be a component in the Bernstein variety $\Omega(G)$, see
\cite{Be1}. The Bernstein decomposition

\[ \mathcal{H}(G) =  \bigoplus \mathcal{H}(\Omega)\]

is a canonical decomposition of $\mathcal{H}(G)$ into $2$-sided
ideals $\mathcal{H}(\Omega)$.    This determines the Bernstein
decomposition of the $C^*$-algebra $A$ (where $A = C^*_r(G)$):

\[ A = \bigoplus A(\Omega).\]

This is a canonical decomposition of $A$ into a $C^*$-direct sum
of $C^*$-ideals $A(\Omega)$.

Let $Irr^t(G)$ denote the tempered dual of $G$.   The Bernstein
decomposition of $A$ determines the Bernstein decomposition of
$Irr^t(G)$:

\[Irr^t(G) = \bigsqcup Irr^t(G)_{\Omega}.\]

It is sometimes convenient to assimilate the component $\Omega$ to
a point $s$, in which case the Bernstein decomposition is

\[A = \bigoplus A(s)\]

where $s$ is a point in the Bernstein spectrum $\mathcal{B}(G)$.

In this article we refine the Bernstein decomposition of $A$. Each
$C^*$-summand $A(s)$ is itself a finite direct sum :

\[A(s) = A(t_1) \oplus \ldots \oplus A(t_j)\]

where $j$ depends on $s$.   This refinement is minimal: the dual
of each $A(t)$ is a compact {\it connected} orbifold.

In \cite{P2} we construct an equivalence bimodule $\mathcal{E}$
which effects a strong Morita equivalence between $A$ and a
commutative $C^*$-algebra $D$. For each $s$, a direct summand
$\mathcal{E}(s)$ of $\mathcal{E}$ effects a strong Morita
equivalence between $A(s)$ and a commutative sub-$C^*$-algebra
$D(s)$ of $D$. These commutative $C^*$-algebras admit a very
simple, uniform description as follows.

Let $s$ be a point in the Bernstein spectrum $\mathcal{B}(G)$. We
can think of $s$ as a vector $(\tau_1, \ldots, \tau_k)$ of
irreducible supercuspidal representations of smaller general
linear groups: the entries of this vector are determined up to
tensoring with unramified quasicharacters and permutation.   If
the vector is $(\sigma_1, \ldots, \sigma_1, \ldots, \sigma_r,
\ldots, \sigma_r)$ with $\sigma_j$ repeated $e_j$ times, $1 \leq j
\leq r$, and $\sigma_1, \ldots, \sigma_r$ pairwise distinct (after
unramified twist) then we say that $s$ has {\it exponents}
$e_1(s), \ldots, e_r(s)$.  Let

\[d(s) = e_1(s) + \ldots + e_r(s)\]

and let

\[W(s) = S_{e_1(s)} \times \ldots \times S_{e_r(s)}\]

a product of symmetric groups.   Let $\mathbb{T}^d$ denote the
standard compact torus of dimension $d$.    We prove the following
strong Morita equivalence:

\[A(s) \sim C(\widehat{\mathbb{T}^{d(s)}}/W(s))\]

the {\it extended} quotient of the torus $\mathbb{T}^{d(s)}$ by
the finite group $W(s)$.

Note that this strong Morita equivalence {\it depends only on the
exponents of $s$}.   If $\Omega_1 \cong \Omega_2$ as complex
algebraic varieties then we will write $s_1 \sim s_2$.   Note that
$\Omega _1 \cong \Omega_2$ if and only if $s_1 \sim s_2$ if and
only if $s_1, s_2$ have the same exponents.   Then $s_1 \sim s_2$
if and only if $A(s_1)$ and $A(s_2)$ are strongly Morita
equivalent:

\[
s_1 \sim s_2 \Longleftrightarrow A(s_1) \sim A(s_2).
\]

The support of Plancherel measure is therefore given by

\[
Irr^t(G) = \bigsqcup_s \widehat{\mathbb{T}^{d(s)}}/W(s).
\]

We do not have a numerical formula for Plancherel density on
$Irr^t(G)$, nor does \cite{BHK}, but we hope that the above
structure theorems will be a useful step on the way towards the
final numerical formula. Meanwhile (in joint work with Paul Baum)
we prove that {\it Plancherel measure is rotation-invariant}. This
implies that Plancherel density is a constant on each circle in
the discrete series of $GL(n)$.

There is one summand in the Bernstein decomposition of $A$ which
is conspicuous, namely the component determined by the cuspidal
pair $(M,1)$ where $M$ is the diagonal subgroup of $G$ and $1$ is
the trivial representation of $M$.   This component is itself
strongly Morita equivalent to the reduced Iwahori-Hecke
$C^*$-algebra $C^*_r(G//I)$.   We prove, in great detail, a
structure theorem for this unital $C^*$-algebra, using standard
techniques of noncommutative Fourier analysis.

The dual of the reduced Iwahori-Hecke $C^*$-algebra $C^*_r(G//I)$
is a compact Hausdorff space which is best viewed as a space of
Deligne-Langlands parameters, as done in \cite{HP}.

Although it is not strictly necessary, the geometry of the
tempered dual $Irr^t(G)$ is best viewed in terms of Langlands
parameters, see \cite{BP2}.

In some ways, the Schwartz algebra $\mathcal{C}(G)$ is more
appropriate to Plancherel measure than the reduced $C^*$-algebra.
This was certainly the point of view of Harish-Chandra \cite{HC},
and also the one adopted in \cite{BHP1,BHP2}. We should also
mention the Harish-Chandra Product Formula for Plancherel measure,
see \cite[p. 92 -- 93]{HC}, \cite[Theorem 5, p.359]{HC}.

This article is an updated version of \cite{P3, P4}, which have
lain dormant for several years. The present article owes much to
noncommutative geometry, specifically the periodic cyclic homology
computations which appear in \cite{BHP1,BHP2,BP1,BP2}.

\section{Bernstein Decomposition}

Let $\Omega$ be a component in the Bernstein variety $\Omega(G)$.
The Hecke algebra $\mathcal{H}(G)$ can be decomposed into a direct
sum of $2$-sided ideals:

\[\mathcal{H}(G) = \bigoplus \mathcal{H}(\Omega).\]

This is the Bernstein decomposition of the Hecke algebra
$\mathcal{H}(G)$, see \cite{Be1}.

Choose left-invariant Haar measure on $G$.   The left regular
representation $\lambda$ of $L^1(G)$ on $L^2(G)$ is defined as
follows: \[ (\lambda(f))(h) = f * \, h,  \] where $f \in L^1(G), h
\in L^2(G)$ and $*$ denotes convolution.   The reduced
$C^*$-algebra is the closure (in the norm topology) of the image
of $\lambda$:
\[
A = C^*_r(G) = \overline{\lambda(L^1(G))} \subset
\mathcal{L}(L^2(G)).
\]
The dual of $A$ is homeomorphic to the tempered dual of $G$, see
\cite{D}.

We have $\mathcal{H}(G) \subset A$. Let $A(\Omega)$ be the norm
closure in $A$ of $\mathcal{H}(G)$. Then $A(\Omega)$ is a
$C^*$-ideal of $A$. Let $\Omega_1, \Omega_2$ be distinct
components of the Bernstein variety $\Omega(G)$, and let $I =
A(\Omega_1), J = A(\Omega_2)$. Now $\widehat{(I \cap J)} = \hat{I}
\cap \hat{J}$ by \cite[3.2.3]{D}. So $I \cap J \neq 0$ if and only
if there exists $\omega \in \hat{I} \cap \hat{J}$ by
\cite[3.2.3]{D}.   But points in the dual of $A(\Omega)$
correspond to tempered representations $\omega$ such that
\[inf.ch.\, \omega \in \Omega.\]

Since the infinitesimal character map is finite-to-one \cite{Be1},
we must have $I \cap J = 0.$

As a consequence, the Bernstein decomposition of $\mathcal{H}(G)$
determines a $C^*$-direct sum of $C^*$-ideals

\[
A = \bigoplus A(\Omega).
\]

This is the Bernstein decomposition of $A$.   It corresponds to
the following disjoint union:

\[
Irr^t(G) = \bigsqcup Irr^t(G)_{\Omega}.
\]

This is a partition of the tempered dual of $G$, the points in
$Irr^t(G)_{\Omega}$ corresponding to tempered representations
$\omega$ such that $inf.ch.\, \omega \in \Omega$.

Let us consider the following special case.   Let $G = GL(n)$, let
$M_0$ be the diagonal subgroup of $G$, and let $1$ be the trivial
representation of $M_0$.   Then $(M_0, 1)$ is a cuspidal pair. Let
$\Omega_0$ be the component containing $(M_0, 1)$.   We shall call
this component the {\it Borel} component.   As an algebraic
variety, the Borel component is the symmetric product of $n$
copies of $\mathbb{C}^{\times}$.    This is a nonsingular complex
affine algebraic variety.

\section{Reduced Iwahori-Hecke $C^*$-algebra for $GL(n)$}

We begin with an elementary account of the {\it extended}
quotient, see \cite{BC}.  Let $X$ be a space on which the finite
group $\Gamma$ acts. The extended quotient associated to this
action is the quotient space $\hat{X}/\Gamma$ where
\[
\hat{X} = \{(\gamma,x) \in \Gamma \times X : \gamma x = x \}.
\]
The group action on $\hat{X}$ is $g.(\gamma, x) = (g\gamma g^{-1},
gx)$.   Let $X^{\gamma} = \{ x \in X : \gamma x = x \}$ and let
$Z(\gamma)$ be the $\Gamma$-centralizer of $\gamma$. Then the
extended quotient is given by:
\[
\hat{X}/\Gamma = \bigsqcup_{\gamma}X^{\gamma}/Z(\gamma)
\]

where one $\gamma$ is chosen in each $\Gamma$-conjugacy class. If
$\gamma = 1$ then $X^{\gamma}/Z(\gamma) = X/\Gamma$ so {\it the
extended quotient always contains the ordinary quotient}:
\[
\hat{X}/\Gamma = X/\Gamma \sqcup \ldots
\]

We shall need only the special case in which $X$ is the compact
torus $\mathbb{T}^n$ of dimension $n$ and $\Gamma$ is the
symmetric group $S_n$ acting on $\mathbb{T}^n$ by permuting
co-ordinates.

Let $\alpha$ be a partition of $n$, and let $\gamma$ have cycle
type $\alpha$.   Each cycle provides us with one circle, and
cycles of equal length provide us with a symmetric product of
circles.   For example, the extended quotient
$\widehat{\mathbb{T}^5}/S_5$ is the following disjoint union of
compact orbifolds (one for each partition of $5$):
\[
\mathbb{T} \bigsqcup \mathbb{T}^2 \bigsqcup \mathbb{T}^2 \bigsqcup
(\mathbb{T} \times Sym^2 \mathbb{T}) \bigsqcup (\mathbb{T} \times
Sym^2 \mathbb{T}) \bigsqcup (\mathbb{T} \times Sym^3 \mathbb{T})
\bigsqcup \mathbb{T}^5
\]
where $Sym^n \mathbb{T}$ is the $n$-fold symmetric product of the
circle $\mathbb{T}$.  We shall see that this extended quotient is
a model of the arithmetically unramified tempered dual of $GL(5)$.

{\bf 3.1} {\sc Theorem}  Let $I$ be the Iwahori subgroup of
$GL(n)$. Let $n = n_1 + \ldots + n_k$ be a partition of $n$, let
$St(n_j)$ be the Steinberg representation of $GL(n_j)$ and let
$\chi_1, \ldots, \chi_k$ be unramified characters of $GL(1)$. Then
the representation
\[
(\chi_1 \circ det) St(n_1) \times \ldots \times (\chi_k \circ det)
St(n_k)
\]

is unitary, irreducible, tempered and admits $I$-fixed vectors.
Moreover, all such representations are accounted for in this way.

\begin{proof}   We use the Langlands
classification for $GL(n)$, see \cite{Ku}.   Each tempered
representation of $GL(n)$ is of the form $Q(\Delta_1) \times
\ldots \times Q(\Delta_k)$ where the Langlands quotient
$Q(\Delta_i)$ is square-integrable for each $i = 1, \ldots, k$. We
now use transitivity of parabolic induction \cite{} and Borel's
theorem \cite{Bo} to infer that (up to unramified unitary twist)
we must have
\[
\Delta_i = \{|\;|_F^{(1 - n_i)/2}, \ldots, | \; |_F^{(n_i -
1)/2}\}
\]

with $i = 1, \ldots, k$.   But then $Q(\Delta_i)$ is the Steinberg
representation $St(n_i)$ of $GL(n_i)$.   Note that $Q(\chi \otimes
\Delta_i) = (\chi \circ det) \otimes Q(\Delta_i)$ and that
$Q(\Delta_1) \times \ldots \times Q(\Delta_k)$ is irreducible.
\end{proof}

{\bf 3.2} Theorem.   The parameter space for the tempered
representations of $GL(n)$ which admit $I$-fixed vectors is the
extended quotient $\widehat{\mathbb{T}^n}/S_n$.

\begin{proof}  Suppose that there are $r_j$ blocks of size $n_j$ with
$1 \leq j \leq l$. Then the Weyl group of the Levi factor $M =
GL(n_1) \times \ldots \times GL(n_l)$ is
\[
W(M) = S_{r_1} \times \ldots \times S_{r_l}.
\]

This Weyl group permutes blocks of the same size.   By standard
Bruhat theory, the Weyl group controls equivalences of
parabolically induced representations.   It follows that the
parameter space for the tempered representations which admit
$I$-fixed vectors is
\[
X = \bigsqcup \mathbb{T}^{r_1}/S_{r_1} \times \ldots \times
\mathbb{T}^{r_l}/S_{r_l}.
\]

The disjoint union is over all partitions
\[
n_1 + \ldots + n_1 + \ldots + n_l + \ldots + n_l = r_1n_1 + \ldots
+ r_ln_l = n.
\]

Let now $\gamma$ be an element in $S_n$ whose cycle type is the
above partition.   Then the centralizer $Z(\gamma)$ is the product
of wreath products:

\[
Z(\gamma) = (\mathbb{Z}/n_1 \wr S_{n_1}) \times \ldots \times
(\mathbb{Z}/n_l \wr S_{n_l}).
\]

But the cyclic groups $\mathbb{Z}/n_1, \ldots, \mathbb{Z}/n_l$ act
trivially on the fixed-point set $(\mathbb{T}^n)^{\gamma}$. We
then have

\begin{eqnarray*}
\widehat{\mathbb{T}^n}/S_n &=&
\bigsqcup(\mathbb{T}^n)^{\gamma}/Z(\gamma)\nonumber\\ &=&
\bigsqcup\{(a, \ldots, a, \ldots , b, \ldots, b, \ldots, c,
\ldots, c, \ldots, d, \ldots,d)\}/W(M)\nonumber\\ &=& \bigsqcup
Sym^{r_1} \mathbb{T} \times \ldots \times Sym^{r_l}
\mathbb{T}\nonumber\\ &=& X.\nonumber\\
\end{eqnarray*}

\end{proof}

Left-invariant Haar measure on $G$ is now chosen so that the
Iwahori subgroup has volume $1$.   Let $e : G \rightarrow
\mathbb{R}$ be defined as follows: $e(x) = 1$ if $x \in I$, $e(x)
= 0$ if $x \notin I$.  Then $e$ is a projection in $A$.   It
generates the non-unital $C^*$-ideal $AeA$.   Now consider the
corner $eAe$.   This is a unital algebra with unit $e$.

{\bf 3.3} {\sc Lemma}  The equivalence bimodule $eA$ effects a
strong Morita equivalence between $eAe$ and $AeA$.

\begin{proof}  Let $B = eAe, C = AeA, \mathcal{E} = eA$.   We have
to check 3 points.

(1)   $\mathcal{E}$ admits a $C$-valued inner product given by
\[
(x,y)_{\mathcal{E}} = x^*y.
\]
Then $||x||_{\mathcal{E}} = ||x^*x||_C^{1/2} = ||x||$.   Since $C
= AeA$ is a $C^*$-ideal, $C$ is $|| \, . \,||$-complete and
$\mathcal{E}$ is a $C^*$-module.

The set $\{(x,y)_{\mathcal{E}} : x, y \in \mathcal{E}\}$ is
\[
\{x^*y : x, y \in \mathcal{E}\} = \{a^*eb : a, b \in A\} = AeA = C
\]
hence $\mathcal{E}$ is a full $C^*$-module.

(2)   $\mathcal{E}$ is a right $C$-module given by $\mathcal{E}
\times C \longrightarrow \mathcal{E}, (ea,c) \mapsto eac$.

(3)   The standard rank $1$ operators are given by
$\theta_{x,y}(z) = x(y,z) = xy^*z$ with $x, y, z \in \mathcal{E}$.
Then $\theta_{ex,e}(z) = (exe)z$ and so the linear span of the
$\theta_{x,y}$ is isomorphic to $eAe$. Since $eAe$ is complete,
the closure of the linear span of the $\theta_{x,y}$ is isomorphic
to $eAe$.  So we have \[ eAe \cong End^0(\mathcal{E})\] the
compact endomorphisms of $\mathcal{E}$. Hence $\mathcal{E}$ is an
equivalence bimodule.
\end{proof}

{\bf 3.4} {\sc Theorem}   The dual of $eAe$ is homeomorphic to the
extended quotient $\widehat{\mathbb{T}^n}/S_n$.
\begin{proof}   The equivalence bimodule
$\mathcal{E}$ determines a homeomorphism of dual spaces:
\[
\widehat{eAe} \cong \widehat{AeA}.
\]
Also we have
\[
\pi(e) = \int e(x)\pi(x) \, dx = \int_I \pi(x) \, dx
\]

which is the projection onto the subspace of $I$-fixed vectors
occurring in $\pi$.   So $\pi(e) \neq 0$ if and only if $\pi$
admits nonzero $I$-fixed vectors.   Since $\pi(exe) =
\pi(e)\pi(x)\pi(e)$, the dual of $eAe$ is precisely that part of
$Irr^t(G)$ which admits nonzero $I$-fixed vectors.   Therefore we
have $AeA = A(\Omega)$ where $\Omega$ is the Borel component in
the Bernstein variety $\Omega(G)$.  We now apply Theorem 3.2.
\end{proof}

Let $\mathcal{H}(G//I)$ be the Iwahori-Hecke algebra.   This
comprises all complex-valued functions $\phi$ on $G$ which are
compactly supported and bi-invariant with respect to $I$:
\[
\phi(i_1xi_2) = \phi(x)
\]
for all $i_1, i_2 \in I, x \in G$.   The product in
$\mathcal{H}(G//I)$ is the convolution product.

Since $\phi = e * \phi * e$ it is immediate that
$\mathcal{H}(G//I) \subset eAe$.   In fact $\mathcal{H}(G//I)$ is
dense in $eAe$ in the reduced $C^*$-algebra norm, and we refer to
$eAe$ as the reduced Iwahori-Hecke $C^*$-algebra.   The notation
is
\[
C^*_r(G//I) = eAe.
\]
A hermitian vector bundle $S$ now presents itself.   The base
space $X$ is the dual of $eAe$.   By Theorem 3.4, the base space
is the extended quotient $\widehat{\mathbb{T}^n}/S_n$.   The total
space $S$ is the set of all $I$-fixed vectors.   The fibre
$S_{\pi}$ comprises all $I$-fixed vectors in the representation
$\pi$.   This bundle is a trivial vector bundle.

The bundle $S$ admits an endomorphism bundle
\[
End(S) = S \otimes S^*.
\]
Since $S$ is hermitian, the continuous sections of $End(S)$ form a
unital $C^*$-algebra whose dual is homeomorphic to $X$.

{\it Definition}.   Let $\phi \in \mathcal{H}(G//I)$. The Fourier
Transform $\hat{\phi}$ is defined as \[ \hat{\phi}(\pi) =
\pi(\phi) = \int \phi(g)\pi(g) dg. \]

Since $\phi = e\phi e$ it follows that
\[
\hat{\phi}(\pi) = \pi(e\phi e) = \pi(e)\pi(\phi)\pi(e)\] where
$\pi(e)$ projects onto the $I$-fixed subspace of $\pi$.  Therefore
we have $\hat{\phi}(\pi) \in End(S_{\pi})$.   Define \[\alpha :
\mathcal{H}(G//I) \longrightarrow C(End \, S)\] \[ \phi \mapsto
\hat{\phi}.\] So we have\[(\alpha(\phi))(\pi) = \hat{\phi}(\pi)
\in End(S_{\pi}).\]

{\bf 3.5} {\sc Theorem}   The Fourier Transform extends uniquely
to an isomorphism of unital $C^*$-algebras:
\[
C^*_r(G//I) \cong C(End \, S).
\]
This is a finite direct sum of homogeneous $C^*$-algebras.

\begin{proof}   We have already shown that the Fourier Transform
determines a map \[ \mathcal{H}(G//I) \longrightarrow C(End \,
S).\] {\it Injectivity}.   Let $y \in eAe$.   Then there exists
$\pi$ in $\widehat{eAe}$ such that $||y|| = sup||\pi(y)||$ by
\cite[3.3.6]{D}.   Therefore $y \neq 0$ implies there exists $\pi$
such that $\pi(y) \neq 0$, and $ y_1 \neq y_2$ implies there
exists $\pi$ such that $\pi(y_1) \neq \pi(y_2)$ so that
\[
\alpha : \mathcal{H}(G//I) \longrightarrow C(End \, S)
\]
is injective.

{\it Surjectivity}.   The image $\alpha(eAe)$ is a
sub-$C^*$-algebra of $C(End \, S)$.   Let $B = \alpha(eAe)$. Note
that $C(End \, S)$ is a liminal $C^*$-algebra with compact
Hausdorff dual $X$.   Let
\[
C = C(End \, S).
\]

We shall now apply \cite[11.1.4]{D}, which is a preliminary
version of the Stone-Weierstrass theorem for $C^*$-algebras.

 Let $\pi \in
\hat{C}$. Then $\pi \in X$. Consider the restriction $\pi|B$. This
is irreducible because
\[
\pi : \mathcal{H}(G//I) \longrightarrow End(S_{\pi})
\]
is a simple $\mathcal{H}(G//I)$-module, $\mathcal{H}(G//I)$ is
dense in $eAe$, and $B$ is isomorphic to $eAe$.  The simplicity of
the $\mathcal{H}(G//I)$-module of $I$-fixed vectors in $\pi$ is a
classical result of Borel \cite{Bo}.

Let $\pi, \psi \in \hat{C}$ and suppose $\pi \neq \psi$.  Then the
simple $\mathcal{H}(G//I)$-modules $S_{\pi}$ and $S_{\psi}$ are
distinct owing to the bijection between irreducible
representations of $G$ which admit nonzero $I$-fixed vectors and
simple $\mathcal{H}(G//I)$-modules. Once again $\mathcal{H}(G//I)$
is dense in $eAe$ which is isomorphic to $B$. So we conclude that
$\pi|B$ and $\psi|B$ are distinct irreducible representations of
$B$.   By \cite[11.1.4]{D}, we have $B = C$.
\end{proof}

There are in fact two equivalence bimodules at work here.  We
recall that the Borel component $\Omega$, as a complex algebraic
variety, is the symmetric product of $n$ copies of
$\mathbb{C}^{\times}$.   The points in the dual of $A(\Omega)$
correspond to tempered representations of $GL(n)$ whose
infinitesimal characters lie in $\Omega$.

The $C^*$-ideal $A(\Omega)$ is given by $A(\Omega) = AeA$. The
first strong Morita equivalence is
\[
AeA \sim eAe
\]
with equivalence bimodule $\mathcal{E} = eA$.   The unital corner
$eAe$ is the reduced Iwahori-Hecke $C^*$-algebra: $eAe =
C^*_r(G//I)$. The second strong Morita equivalence is
\[
eAe \sim C(X)
\]
where $X$ is the extended quotient $\widehat{\mathbb{T}^n}/S_n$.
The equivalence bimodule comprises all continuous sections of the
complex Hermitian bundle $S$ of all Iwahori fixed vectors.

So we have
\[
 A(\Omega) \sim eAe \sim C(X).
\]

The extended quotient $X$ is a disjoint union: one connected
component is the ordinary quotient $\mathbb{T}^n/S_n$.   The
corresponding direct summand of $A(\Omega)$ is discussed in the
next section.

\section{Reduced spherical $C^*$-algebra for $GL(n)$}

Let $K$ be a maximal compact subgroup of $GL(n)$.   We may take $K
= GL(n, \mathcal{O})$.   We have
\[
C^*_r(G//K) \subset C^*_r(G//I).
\]
In the strong Morita equivalence
\[
C^*_r(G//I) \sim C(X)
\]
the sub-$C^*$-algebra $C^*_r(G//K)$ determines a strong Morita
equivalence
\[
C^*_r(G//K) \sim C(\mathbb{T}^n/S_n)
\]
with equivalence bimodule $C(L)$, where $L$ is the line bundle of
all $K$-fixed vectors.   This of course means that
\[
C^*_r(G//K) \cong C(\mathbb{T}^n/S_n).
\]
This is the $C^*$-algebra version of the Satake isomorphism
\cite[p. 147]{Ca}:

\[
\mathcal{H}(G//K) \cong \mathbb{C}[\Lambda]^W
\]

where the lattice $\Lambda = \mathbb{Z}^n$.   The group algebra
$\mathbb{C}[\Lambda]$ will Fourier Transform to a dense subalgebra
of $C(\mathbb{T}^n)$ and the $W$-invariant part
$\mathbb{C}[\Lambda]^W$ will Fourier Transform to a dense
subalgebra of $C(\mathbb{T}^n/W)$.

The $C^*$-algebra $C^*_r(G//K)$ is the reduced {\it spherical}
$C^*$-algebra for the group $GL(n)$, cf. \cite[4.4]{Ca}.   We
conclude with the following

{\bf 4.1} {\sc Theorem}   Let $C^*_r(G//K)$ be the reduced
spherical $C^*$-algebra for $GL(n)$.   Then we have
\[
C^*_r(G//K) \cong C(\mathbb{T}^n/S_n).
\]

\section{Structure of the $C^*$-summand $A(s)$}

Let $s$ be a point in the Bernstein spectrum of $GL(n)$ and let
$A(s)$ be the corresponding  $C^*$-summand in the Bernstein
decomposition of $A$.

{\bf 5.1} {\sc Theorem}  The $C^*$-algebra $A(s)$ is strongly
Morita equivalent to a commutative $C^*$-algebra:
\[A(s) \sim C(\widehat{\mathbb{T}^{d(s)}}/W(s)).\]

\begin {proof}   We begin by reviewing the Bernstein variety
$\Omega(G)$, see \cite{Be1}.   The Bernstein variety is a disjoint
union of ordinary quotients:

\[
\Omega(G) = \bigsqcup D/W(M,D)
\]

where $D$ is a complex torus, and $W(M,D)$ is a certain finite
group which acts on $D$.   We now replace, as in \cite{BP2}, the
ordinary quotient by the extended quotient to create a new variety
$\Omega^+(G)$.   So we have

\[
\Omega^+(G) = \bigsqcup \hat{D}/W(M,D).
\]

We call $\Omega^+(G)$ the {\it extended} Bernstein variety.  Each
component in the extended Bernstein variety is the quotient of a
complex torus by a finite group, and so is itself a complex
algebraic variety.   So $\Omega^+(G)$ is a complex algebraic
variety with countably many irreducible components.

In \cite{BP2} we construct a bijection:

\[
Irr(GL(n) \longrightarrow \Omega^+(GL(n))
\]

This bijection is constructed in terms of Langlands parameters,
and also depends on a short but intricate piece of combinatorics
\cite[p. 217]{BP2}.    Then, by transport of structure, the smooth
dual acquires the structure of complex algebraic variety with
countably many irreducible components.

We will write

\[
\Omega = D/W(M,D)
\]

for a component in $\Omega(G)$.   If we now assimilate $\Omega$ to
a point $s$ in the Bernstein spectrum,  then we can write
\[
\hat{D}/W(M,D) = \widehat{T^{d(s)}_c}/W(s)
\]

where $T^n_c$ denotes the complex torus of dimension $n$.

Now the complex commutative Lie group $T^n_c$ maps by a
deformation retraction onto its maximal compact subgroup
$\mathbb{T}^n$:

\[
(z_1, \ldots, z_n) \mapsto (|z_1|^{-1}z_1, \ldots, |z_n|^{-1}z_n).
\]

There is also a deformation retraction of $Irr(G)$ onto the
tempered dual $Irr^t(G)$, see \cite[Theorem 2]{BP2}.    These
deformation retractions are compatible in the sense that the
following diagram is commutative:

\[
\begin{CD}
Irr(G)  @>>>    \bigsqcup \widehat{T^{d(s)}_c}/W(s)\\ @VVV @VVV\\
Irr^t(G)   @>>>  \bigsqcup \widehat{\mathbb{T}^{d(s)}}/W(s)\\
\end{CD}
\]

\bigskip

In this diagram, the horizontal maps are bijective.

We now return to $C^*$-algebras. The main result of \cite{P1} is
the following isomorphism of $C^*$-algebras:

\[
C^*_r(G) \cong C_0(Irr^t(G), \mathcal{K}(H))
\]

where $\mathcal{K}(H)$ is the $C^*$-algebra of compact operators
on the standard Hilbert space $H$.   For each point $s$ in the
Bernstein spectrum we then have

\[
A(s) \cong C(Irr^t(G)_{\Omega}, \mathcal{K}(H)).
\]

The equivalence bimodule $\mathcal{E} = C(Irr^t_{\Omega}(G), H)$
then effects the following strong Morita equivalence:

\[
A(s) \sim C(\widehat{\mathbb{T}^{d(s)}}/W(s)).
\]
\end{proof}

Let $\mathcal{C}(G)$ denote the Schwartz algebra of $G$.  Let $S$
be a complete set of standard tori in $G$ no two of which are
conjugate in $G$.  The following decomposition is due to
Harish-Chandra \cite[p. 367]{HC}:

\[
\mathcal{C}(G) = \bigoplus_{A \in S} \mathcal{C}_A(G).
\]

There is a refinement of this decomposition:

\[
\mathcal{C}(G) = \bigoplus \mathcal{C}(\mathfrak{o})
\]

where $\mathfrak{o}$ is an $\mathfrak{a}^*$-orbit in $E_2(M)$, see
\cite[p. 359]{HC}. This "wave-packet" decomposition is not stated
explicitly by Harish-Chandra, but is implied by his formula for
the component $f_{\mathfrak{o}}$, see \cite[p. 360]{HC}.

\section{Rotation invariance of Plancherel measure}

This section is joint work with Paul Baum.   Let $G$ be a
reductive $p$-adic group and let $\Psi(G)$ denote the group of
unramified quasicharacters of $G$. It acts naturally on $Irr(G)$
by $\psi : \pi \mapsto \psi \otimes \pi$. As
 in \cite{BDK}, define an algebraic action of $\Psi(G)$ on
$\Omega(G)$ by $\psi : (M, \rho) \mapsto (M, \psi |_M \, . \,
\rho)$.   Then $inf.\, ch$ is a $\Psi(G)$-equivariant map.

If $M$ is a standard Levi subgroup, then $\Psi(M)$ is a complex
commutative Lie group, and its maximal compact subgroup is denoted
$\Psi^t(M)$, the group of unramified (unitary) characters of $M$.

We recall the definition of the Harish-Chandra parameter space
$\Omega^t(G)$. We call a {\it discrete-series pair} a pair $(M,
\sigma)$ where $M$ is a standard Levi subgroup, $\sigma$ is an
irreducible unitary representation in the discrete series of $M$.
We denote by $\Omega^t(G)$ the set of all discrete-series pairs up
to conjugation by $G$. For any discrete-series pair $(M, \sigma)$
the image of the map $\Psi^t(M) \rightarrow \Omega^t(G)$, given by
$\chi \mapsto (M, \chi\sigma)$, is called a {\it connected
component} of $\Omega^t(G)$. Each connected component is a compact
connected orbifold.

The group $\Psi^t(G)$ acts naturally on $Irr^t(G)$ by $\chi : \pi
\mapsto \chi \otimes \pi$.  We now define a smooth action of
$\Psi^t(G)$ on $\Omega^t(G)$ by $\chi : (M, \sigma) \mapsto (M,
\chi |_M \, . \, \sigma)$.

We emphasize that
\[
inf.\, ch : Irr^t(G) \longrightarrow \Omega^t(G) \] is a
homeomorphism of the tempered dual in its standard topology onto
the Harish-Chandra parameter space in its natural topology.  Also,
$inf. \, ch$ is a $\Psi^t(G)$-equivariant map.

Now $\Omega^t(G)$ is the support of Plancherel measure, and
$\Psi^t(G)$ is a compact torus.   We shall say that $\Psi^t(G)$
{\it rotates} $\Omega^t(G)$, and that it {\it rotates} Plancherel
measure $\nu$.   If $\chi \in \Psi^t(G)$ then the rotated measure
is given by
\[
E \mapsto \nu(\chi^{-1}(E)
\]

where $E$ is a Borel set in the tempered dual of $G$.

If $G = GL(n)$, then $\Psi^t(G) = \mathbb{T}$.   The
representations in the discrete series of $G$ arrange themselves
into circles.

{\bf 6.1} {\sc Theorem}.   Let $G$ be a reductive $p$-adic group.
Then Plancherel measure is rotation-invariant.  In the case of
$GL(n)$, Plancherel measure induces Haar measure on each circle in
the discrete series.

\begin{proof}

We begin with the Plancherel formula

\[
f(1) = \int tr \, \hat{f}(\omega) \, d\nu(\omega)
\]

where $d\nu$ is Plancherel measure.   We have

\begin{eqnarray}
\hat{f}(\chi.\omega) &=& \int f(g)(\chi.\omega)(g) \,dg
\nonumber\\ &=& \int f(g) \chi(g) \omega(g) \, dg \nonumber\\ &=&
\widehat{f.\chi}(\omega).\label{omega}
\end{eqnarray}

Let $f \in \mathcal{C}(G)$.   Making a change of variable, and
applying (\ref{omega}) we get

\begin{eqnarray*}
\int tr \, \hat{f}(\omega) d\nu (\chi^{-1}.\,
\omega) &=& \int tr
\, \hat{f}(\chi.\omega) d\nu(\omega)\nonumber\\ &=& \int tr \,
\widehat{f.\chi}(\omega) d\nu(\omega)\nonumber\\ &=&
(f.\chi)(1)\nonumber\\ &=& f(1).\nonumber\\
\end{eqnarray*}

We now apply the trace-Paley-Wiener theorem \cite{BDK} and
properties of the $q$-projection (twisted projection) in
\cite{BP2}.   As $f$ varies in the Hecke algebra $\mathcal{H}(G)$
the integrand $tr \, \hat{f}(\omega)$ varies over all finite
Fourier series on the component $\Omega^t(G)$.   So the Radon
measures $d\nu(\omega)$ and $d\nu(\chi^{-1}.\omega)$ determine the
same continuous linear functional on a dense subset of
$C_0(\Omega^t(G))$.   Therefore these two Radon measures are
equal, i.e. $\nu$ is rotation-invariant.

In the case of $GL(n)$, each circle in the discrete series is
rotated by the group $\Psi^t(G)$, hence Plancherel measure induces
Haar measure on each circle.

\end{proof}

Let $\Omega^t$ be a component in $\Omega^t(G)$.   Then $\Omega^t =
T/W(M,T)$.   Let $d$ be the depth of $M$, i.e. the dimension of
$\Psi(M)$.   Then $T = \mathbb{T}^d$.   When $G = GL(n)$ we have
$\Psi^t(G) \cong \mathbb{T}$ and the action of $\Psi^t(G)$ on
$\Omega^t$ is induced by the {\it diagonal} action of $\mathbb{T}$
on $\mathbb{T}^d$.   So Plancherel density is invariant under this
diagonal action.

\end{document}